# DIMENSION REDUCTION BASED ON CONSTRAINED CANONICAL CORRELATION AND VARIABLE FILTERING[1]

By Jianhui Zhou and Xuming He

*University of Virginia and University of Illinois*

The "curse of dimensionality" has remained a challenge for high-dimensional data analysis in statistics. The sliced inverse regression (SIR) and canonical correlation (CANCOR) methods aim to reduce the dimensionality of data by replacing the explanatory variables with a small number of composite directions without losing much information. However, the estimated composite directions generally involve all of the variables, making their interpretation difficult. To simplify the direction estimates, Ni, Cook and Tsai [*Biometrika* **92** (2005) 242–247] proposed the shrinkage sliced inverse regression (SSIR) based on SIR. In this paper, we propose the constrained canonical correlation ($C^3$) method based on CANCOR, followed by a simple variable filtering method. As a result, each composite direction consists of a subset of the variables for interpretability as well as predictive power. The proposed method aims to identify simple structures without sacrificing the desirable properties of the unconstrained CANCOR estimates. The simulation studies demonstrate the performance advantage of the proposed $C^3$ method over the SSIR method. We also use the proposed method in two examples for illustration.

**1. Introduction.** For data sets with a large number of variables, the well-known "curse of dimensionality" poses a challenge to most statistical methods. Dimension reduction methods are often used to reduce dimensionality, enabling regression or classification to be performed in a parsimonious way.

We consider a regression setting where the univariate response $y$ is related to the explanatory variable $x_{p\times 1}$ through $K$ linear combinations of $x$ and an unknown function $f$. One way to describe the model is

$$y = f(x^T\beta_1, x^T\beta_2, \ldots, x^T\beta_K, \epsilon), \tag{1}$$

Received May 2007; revised July 2007.
[1]Supported in part by NSF Grant DMS-06-04229.

where $\epsilon$ is the random error independent of $x$, $K$ is the smallest integer for model (1) to hold, and $\beta_i$ are a set of effective dimension reduction (e.d.r.) directions as in [11]. Following [2], the subspace $S_{y|x} = \text{span}\{\beta_1, \ldots, \beta_K\} \subseteq \mathbb{R}^p$ is called the central dimension reduction subspace (DRS). The central DRS is unique, even though the sets of e.d.r. directions, which span the central DRS, can be taken differently when $K \geq 2$. Given the central DRS $S_{y|x}$, we use the $K$-dimensional projection of $x$ onto $S_{y|x}$, instead of $x$ itself, in the model.

To estimate the dimensionality $K$ and the central DRS, the sliced inverse regression (SIR) method is proposed by [11] to summarize the explanatory variables into a smaller number of linear projections. Compared to the principal component analysis that summarizes the information contained in the explanatory variables, SIR takes into account extra information contained in the response variable. The sliced average variance estimation (SAVE) method is proposed by [5], and it is shown that it captures a larger portion of the central DRS than SIR by [4]. To achieve exhaustive estimation of the central DRS, simple contour regression (SCR) and general contour regression (GCR) are proposed by [9]. Assuming an additive random error in model (1), the (conditional) minimum average variance estimation (MAVE) method is proposed in [20] by minimizing an objective function that involves both the direction estimation and the nonparametric function estimation. Other dimension reduction methods in the literature include the ordinary least squares (OLS) by [14], the principal Hessian directions (PHD) by [12], and the canonical correlation (CANCOR) method by [8].

When $K < p$, a reduction in dimensionality is achieved through model (1), but each effective direction usually involves all of the explanatory variables. When $p$ is large and the variables in $x$ are of different scales, the estimated linear combinations $x^T \hat{\beta}_i$ are difficult to interpret, which makes them less useful in the analyses following dimension reduction. Attempts have been made to address this problem. For single-index models, which are a special case of model (1), an AIC-based criterion is proposed by [16] to select the relevant variables and the smoothing parameter for estimating the unknown function simultaneously. In the general framework of model (1), the shrinkage sliced inverse regression (SSIR) method is developed in [17] by employing the LASSO approach of [19]. In SSIR, the central DRS $S_{y|x}$ is estimated by $\text{span}\{\text{diag}(\tilde{\alpha})\hat{B}\}$, where $\text{span}\{M\}$ denotes the subspace spanned by the columns of the matrix $M$, $\hat{B}$ corresponds to the estimated central DRS $\text{span}\{\hat{B}\}$ of SIR in [3], and the shrinkage indices $\tilde{\alpha}$ are determined through a LASSO regression subject to an $L_1$-norm constraint on $\tilde{\alpha}$.

In this paper, we propose an approach to reduce the number of variables appearing in the CANCOR directions, aiming to explore possible sparsity in the e.d.r. directions through both dimension reduction and variable selection. We describe the constrained canonical correlation ($C^3$) method in Section 2 by imposing the $L_1$-norm constraint on the e.d.r. direction estimates in CANCOR. In Section 3, a variable filtering procedure is proposed to threshold the estimated coefficients to reduce the number of nonzero coefficients in the direction estimates by $C^3$. The final values of the nonzero coefficients are then re-estimated by CANCOR on the relevant subset of variables. Simulation studies are conducted in Section 4, and two data sets, the car price data and the Boston housing data, are analyzed in Section 5. Concluding remarks about the proposed method can be found in Section 6.

## 2. Constrained canonical correlation.

2.1. *CANCOR revisited.* Using the B-spline basis functions generated for the response variable, the CANCOR method, which is asymptotically equivalent to SIR, is proposed by [8]. Suppose that the range of $y$ is a bounded interval $[a, b]$. Given $k_n$ internal knots in $[a, b]$ and the spline order $m$, we generate $m + k_n$ B-spline basis functions. Under the linearity condition of [11], CANCOR estimates a set of e.d.r. directions by estimating the canonical variates between the B-spline basis functions and $x$. Since the generated $m+k_n$ B-spline basis functions sum to 1, we use in CANCOR the first $m+k_n-1$ basis functions of $y$, $\pi(y) = (\pi_1(y), \ldots, \pi_{m+k_n-1}(y))^T$. Given $n$ observations, $(Y_i, X_i)$, of the random variables $(y, x)$, let $X_{n\times p} = (X_1, \ldots, X_n)^T$ and $\Pi_{n\times(m+k_n-1)} = (\pi(Y_1), \ldots, \pi(Y_n))^T$ be the two data matrices containing the predictor values and the B-spline basis function values. The CANCOR method is then to estimate the canonical correlations between the columns of $X$ and the columns of $\Pi$. The dimensionality of the central DRS is selected by performing the following sequential tests on the number of nonzero canonical correlations, $H_{0,s}: \gamma_s > \gamma_{s+1} = 0$ versus $H_{1,s}: \gamma_{s+1} > 0$ for $s = 0, 1, \ldots, p-1$, where $\gamma_s$ are the asymptotic canonical correlations between $\pi(y)$ and $x$ in decreasing order. For details on the test statistic, see [8]. The dimensionality estimate $\hat{K}$ is the smallest $s$ such that $H_{0,s}$ is not rejected. The directions $\hat{\beta}_i$ in the estimated canonical variates $x^T \hat{\beta}_i$, corresponding to the nonzero correlations, are the estimated e.d.r. directions. The estimated central DRS is $\hat{S}_{y|x} = \text{span}\{\hat{\beta}_1, \ldots, \hat{\beta}_{\hat{K}}\}$.

Using the standardized variables $z = \Sigma_{xx}^{-1/2}[x - E(x)]$, where $E(x)$ and $\Sigma_{xx}$ are mean and covariance of $x$, it is shown in [8] that CANCOR is based on $\text{Cov}[E(z|y)]$, which is the same as the kernel matrix in SIR. Letting $(\lambda_i, \eta_i)$, $i = 1, 2, \ldots, p$, be the eigenvalues in decreasing order and the corresponding eigenvectors of the matrix $\text{Cov}[E(z|y)]$, the canonical correlations between $\pi_i(y)$ and $x$ are $\gamma_i = \lambda_i^{1/2}$, and the canonical directions $\beta_i = \Sigma_{xx}^{-1/2}\eta_i$, corresponding to the nonzero eigenvalues $\lambda_i$, are contained in the central DRS $S_{y|x}$, assuming the linearity condition of [11].

The CANCOR method actually solves an optimization problem that sequentially finds the directions $\beta_i$ with the maximum correlations between $x^T\beta_i$ and some functions of $y$. If the canonical correlations are strictly monotone, then the CANCOR approach leads to the *unique* identification of the directions $\beta_i$, not just the space spanned by those directions. The directions $\beta_i$ have their own interpretations, regardless of whether model (1) holds.

2.2. *Constrained CANCOR.* To estimate the first constrained e.d.r. direction, we solve the following nonlinear constrained canonical correlation ($C^3$) problem, where $\|\beta\|_{L_1}$ will be used throughout the paper for the $L_1$-norm of any vector $\beta$.

PROBLEM 1.

$$\text{Maximize} \quad \alpha^T \hat{\Sigma}_{\pi x}\beta,$$
$$\text{subject to} \quad \alpha^T\hat{\Sigma}_{\pi\pi}\alpha = 1; \qquad \beta^T\hat{\Sigma}_{xx}\beta = 1; \qquad \|\beta\|_{L_1} \leq t.$$

In Problem 1, the matrix $\hat{\Sigma}_{\pi x}$ is taken to be the sample covariance matrix between $\pi(y)$ and $x$, and $\hat{\Sigma}_{\pi\pi}$ and $\hat{\Sigma}_{xx}$ the sample variance-covariance matrices of $\pi(y)$ and $x$, respectively. The maximizer of Problem 1 for a given $t$ is denoted as $(\hat{\alpha}_1^c, \hat{\beta}_1^c)$, and the maximal correlation is $\hat{\gamma}_1^c = \hat{\alpha}_1^{cT}\hat{\Sigma}_{\pi x}\hat{\beta}_1^c$. The criterion of selecting the value of $t$ will be discussed later in this section.

One salient point about Problem 1 is that it does not use the inverse of the sample covariance matrix $\hat{\Sigma}_{xx}$. Even if the number of variables $p$ is greater than $n$, the constrained problem will generally have a unique solution for small $t$, while the unconstrained problem admits possibly infinitely many solutions. When there are multiple candidates to achieve the same correlation, the constrained optimization favors the candidate directions with sparse coefficients.

For the identifiability of the $i$th ($i \geq 2$) constrained e.d.r. directions, we require them to be uncorrelated with the previously estimated directions, that is, the $i$th ($i \geq 2$) constrained e.d.r. direction is estimated by solving the following constrained optimization problem for $(\hat{\alpha}_i^c, \hat{\beta}_i^c)$, and $\hat{\gamma}_i^c = \hat{\alpha}_i^{cT}\hat{\Sigma}_{\pi x}\hat{\beta}_i^c$.

PROBLEM 2.

$$\text{Maximize} \quad \alpha^T\hat{\Sigma}_{\pi x}\beta,$$
$$\text{subject to} \quad \alpha^T\hat{\Sigma}_{\pi\pi}\alpha = 1; \qquad \beta^T\hat{\Sigma}_{xx}\beta = 1; \qquad \|\beta\|_{L_1} \leq t;$$
$$\alpha^T\hat{\Sigma}_{\pi\pi}\hat{\alpha}_l^c = 0; \qquad \beta^T\hat{\Sigma}_{xx}\hat{\beta}_l^c = 0; \qquad l = 1,\ldots,i-1.$$

Starting from Problem 1 and solving Problem 2 iteratively over $i$, we get a set of the constrained e.d.r. direction estimates $\{\hat{\beta}_i^c\}$ for $i = 1, \ldots, \hat{K}$, where $\hat{K}$ is the number of significant e.d.r. directions determined by the sequential tests in CANCOR or any pre-determined dimensionality for the problem. The tuning constant $t$ is expected to vary with $i$, as detailed in the next subsection.

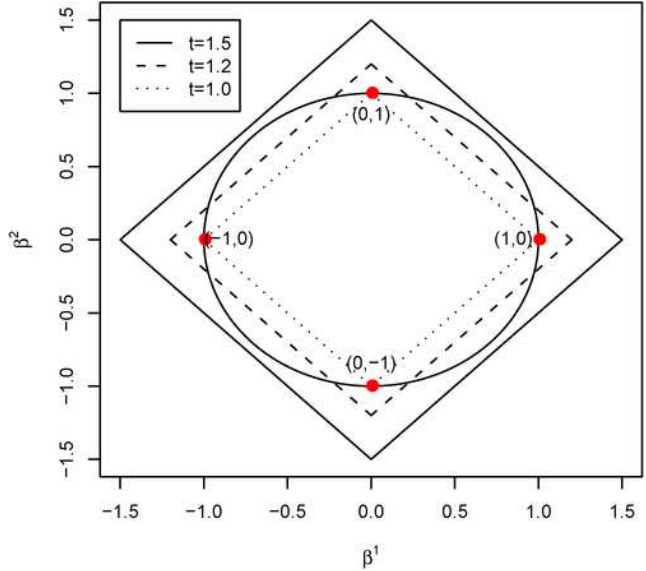

Fig. 1. *Constraint effect in 2D space.*

Note that without the $L_1$-norm constraint $\|\beta\|_{L_1} \leq t$, Problems 1 and 2 indeed estimate the canonical correlations between $\pi(y)$ and $x$. Under the constraint, some coefficients in the direction estimates $\hat{\beta}_i^c$ shrink toward 0. To illustrate that the $L_1$-norm constraint favors sparse coefficients of $\hat{\beta}_i^c$, we discuss the effect of the constraint in Problem 1 in the 2-dimensional space as shown in Figure 1.

In Figure 1, we assume $\hat{\Sigma}_{xx} = I_2$ and $\beta = (\beta^1, \beta^2)$. The circle is the contour defined by $\beta^T \hat{\Sigma}_{xx} \beta = (\beta^1)^2 + (\beta^2)^2 = 1$, and the squares are defined by $|\beta^1| + |\beta^2| = t$ for $t = 1.5$, 1.2 and 1.0, respectively. In CANCOR, we search for the estimate of $\beta$ on the whole circle to maximize the correlation. In $C^3$, the feasible region of $\beta$ is the part of the circle lying in the square defined by $t$. When $t$ gets smaller, the coefficients of the points in the feasible region are closer to being sparse. As shown in Figure 1, when $t$ is decreased to 1, the feasible region shrinks to the four points $(\pm 1, 0)$ and $(0, \pm 1)$. Similarly in the general $p$-dimensional space, when $t$ gets smaller, the points in the feasible region defined by $t$ move toward sparser coordinates.

For selecting the tuning parameter $t$ in Problem 1 or 2, we first study the range of $t$ in which the $L_1$-norm constraint takes effect. When each explanatory variable in $x \in \mathbb{R}^p$ is standardized to variance 1, the matrix $\hat{\Sigma}_{xx}$ has diagonal elements 1 and off-diagonal elements between $-1$ and 1. Thus, we have $\|\beta\|_{L_1}^2 \geq \beta^T \hat{\Sigma}_{xx} \beta = 1$. When $t = 1$, only the unit vectors, such as $(0, \ldots, 0, \pm 1, 0, \ldots, 0)$, are in the feasible region. When $t < 1$, the feasible region is empty and the optimization problem does not have a solution. Therefore, the tuning parameter $t$ should not be less than 1 for a meaningful $L_1$-norm constraint. On the other hand, when $t \geq t_0$, where $t_0 = \|\hat{\beta}_i\|_{L_1}$ and $\hat{\beta}_i$ is the unconstrained e.d.r. direction estimate by CANCOR, the maximizer is always $\hat{\beta}_i$. Therefore, the tuning parameter $t$ should be in the interval $[1, t_0]$.

2.3. *Computational issues.* Starting from $t = t_0$, we search for the value of $t$ in each problem iteratively by decreasing $t$ by a small amount $\Delta t$ in each iteration. Given the value of $t$ in each iteration, the constrained direction $\hat{\beta}_i^c$ and the maximum correlation $\hat{\gamma}_i^c$ are estimated. When $t$ gets smaller, more coefficients in $\hat{\beta}_i^c$ tend to shrink toward 0. As $t$ is decreased to 1, only one variable will be involved in the estimated $x^T \hat{\beta}_i^c$. However, we usually need to stop decreasing $t$ before it reaches 1, because the maximum correlation may have decreased too much when $t$ is too small. To balance the simplicity of the estimated direction with the loss in the correlation, we propose to stop decreasing $t$ when the maximum correlation drops below a lower confidence limit of the canonical correlation. The limiting distributions of the sample canonical correlations can be found in [15]. However, those distributions do not apply directly in our setting. The positive news is that the accuracy (in confidence level) of such a confidence limit is not so critical for our purpose, and we prefer to use easy approximations based on the following Fisher's transformation. Given the estimated canonical correlation $\hat{\gamma}_i$, we transfer it to $\hat{\rho}_i = (1/2) \log[(1 + \hat{\gamma}_i)/(1 - \hat{\gamma}_i)]$, whose distribution could be approximated by the normal distribution with standard deviation $1/\sqrt{n-3}$. Thus, we use $\hat{\rho}_i - Z_{1-\alpha}/\sqrt{n-3}$ as the $100(1-\alpha)\%$ lower confidence limit of the transformed canonical correlation, where $Z_{1-\alpha}$ is the $100(1-\alpha)\%$ quantile of the standard normal distribution. Accordingly, the lower confidence limit of the $i$th canonical correlation could be approximated by $(\exp(2\tau_i) - 1)/(\exp(2\tau_i) + 1)$, where $\tau_i = \hat{\rho}_i - Z_{1-\alpha}/\sqrt{n-3}$. We stop decreasing $t$ in the iterative process when the corresponding maximum correlation falls below the lower confidence limit given above. The value of $t$ in the iteration just before the above process is stopped is selected for Problem 1 or 2, and the corresponding direction estimate $\hat{\beta}_i^c$ is taken as the estimated $i$th constrained e.d.r. direction. The two parameters, $\Delta t$ and $\alpha$, are user-specified. In this paper, we use as default $\alpha = 0.005$ and $\Delta t = 0.05$.

Problems 1 and 2 are nonlinear constrained optimization problems with an $L_1$-norm constraint. Most of the existing algorithms for solving the constrained optimization problems need a pair of initial directions for $\alpha$ and $\beta$, and there is no guarantee that the global maximum will be found through iteration. In the $C^3$ method, we use the unconstrained e.d.r. directions $\hat{\alpha}_i$ and $\hat{\beta}_i$ estimated by CANCOR as the initial directions for solving Problems 1 and 2. In the iterative process for selecting the value of $t$, the constrained direction estimates in the previous iteration are used as the initial directions in the next iteration when $t$ is deceased by $\Delta t$. Since the unconstrained e.d.r. directions correspond to the global maxima, the proposed scheme greatly reduces the chance of hitting a local maximum in the iterations. Furthermore, we always know what to expect from the maximum correlation, so it is easy to know when a poor local maximum is found.

2.4. *Consistency.* Throughout we assume that $\Sigma_{xx} > 0$, and that the nonzero eigenvalues of $\text{Cov}[E(z|y)]$ are distinct, that is, $\gamma_1 > \gamma_2 > \cdots > \gamma_K > 0$. In this setting, the canonical directions $\beta_i$ ($i \leq K$) are uniquely identified, in contrast to those in model (1). Obviously, the directions $\beta_i$ defined through CANCOR are not necessarily part of those in model (1), unless a well-known linearity condition holds on the distribution of $x$ (see [11]). We believe that the composite directions $\beta_i$ are useful even when model (1) does not hold. Letting $k_n$ be the number of internal knots used in the spline, we have:

THEOREM 1. *Assume that $\Sigma_{xx} > 0$, $\gamma_1 > \cdots > \gamma_K > 0$ for some $K$, and that the following conditions $A_1$–$A_4$ hold.*
*$A_1$: The marginal density of $y$ is bounded away from 0 and infinity on $[a,b]$.*
*$A_2$: $k_n \to \infty$ and $k_n = o(n^{2/3})$ as $n \to \infty$.*
*$A_3$: $E(\|x\|^4) < \infty$.*
*$A_4$: Each component of $E(z|y)$ is a function on $[a,b]$ with bounded derivative;*
*Then we have $\hat{\beta}_i^c \xrightarrow{P} \beta_i$ for $i \leq K$.*

Note that conditions $A_1$–$A_4$ are sufficient for the consistency of the constrained canonical correlation estimates, but the assumption of distinct eigenvalues $\gamma_i$, $i = 1, \ldots, K$ is there to ensure the consistency of the constrained direction estimates. Since the direction estimates from $C^3$ will not be taken as the final estimates, we shall defer the results on asymptotic distributions to the next section.

**3. Variable filtering and final estimates.** As demonstrated in Figure 1, the $L_1$-norm constraint helps to shrink the e.d.r. direction estimates toward the sparse regions in $\mathbb{R}^p$, but exact zero coefficients in $\hat{\beta}_i^c$ are not aimed at. To remove the variables with no or very marginal effects on the e.d.r. directions, we propose to threshold the coefficients in $\hat{\beta}_i^c$ through a variable filtering procedure.

The proposed variable filtering procedure is simple and iterative. At each iteration, one more coefficient of $\hat{\beta}_i^c$ is set to be 0 according to the magnitudes of the coefficients in absolute value. Recall that the $x$ variables are standardized individually. Given the constrained e.d.r. direction estimates $\hat{\alpha}_i^c$ and $\hat{\beta}_i^c$ from Problem 1 or Problem 2, the proposed variable filtering procedure involves the following steps:

1. Let $d = p$.
2. Define a new direction $\hat{\beta}_i'^c(d)$ by keeping the largest $d$ coefficients of $\hat{\beta}_i^c$ in absolute value and setting the other $(p - d)$ coefficients to be 0. Find $\hat{\beta}_i^p(d)$ as the projection of $\hat{\beta}_i'^c(d)$ into the space $B_i$, the set of all $\beta$ such that: (i) $\beta^T \hat{\Sigma}_{xx} \hat{\beta}_j^c = 0$, $j \leq i - 1$, (ii) the set of zero coefficients in $\beta$ is the same as that in $\hat{\beta}_i'^c(d)$ and (iii) $\beta^T \hat{\Sigma}_{xx} \beta = 1$.
3. Compute the correlation, $r_d = \text{Corr}(\pi(y)^T \hat{\alpha}_i^c, x^T \hat{\beta}_i^p(d))$, and the BIC-type criterion, $BIC(d) = n \log(1 - r_d^2) + d \log(n)$.
4. Let $d = d - 1$. Repeat steps 2–4 until $d = 0$.

Performing the above variable filtering procedure, we get a sequence of $BIC(d)$ as $d$ decreases from $p$ to 0. Let $d_0$ be the integer at which $BIC(d)$ is minimized. Then, the $p - d_0$ smallest coefficients of $\hat{\beta}_i^c$ in absolute value are set to 0. This proposed variable filtering procedure is a simplified variable selection procedure. In variable filtering, at most $p$ possibilities are considered, which makes it feasible to do even when $p$ is large.

Whenever $d < i$ in Step 2, the projection $\hat{\beta}_i^p(d)$ is usually a zero vector, and in those cases we set $r_d = 0$. Therefore, the number of nonzero coefficients for $\beta_i$ will be no fewer than $i$ by our variable filtering algorithm. This is not an undesirable feature, as it might be unwise to be aggressively dropping variables if a large number of dimensions will be needed in the model. In most applications, we may be looking for the cases of $K \ll p$; the restriction of $d \geq i$ on the results for $\beta_i$ is mild.

It is important to note that the proposed variable filtering procedure should be applied to the constrained direction estimates by the $C^3$ method. If it is applied directly to the unconstrained direction estimates, it will not be effective, as shown in the first example in Section 5, mainly because the sizes of the coefficients are often "distorted" by collinearity of the variables in the unconstrained direction estimates. Variable selection in the unconstrained problem would be a harder problem.

To aim for higher correlation, we propose re-estimation of the directions after variable filtering: the variables selected for any of the first $\hat{K}$ constrained directions are combined, and the CANCOR method is performed on those variables alone. Then, the final estimate of $\beta_i$, now denoted as $\hat{\beta}_i^f$, is obtained by taking the estimates from re-estimation for its nonzero coefficients and by keeping the zero coefficients as determined in the variable filtering stage. The following theorem concerns the asymptotic behavior of the final direction estimates.

THEOREM 2. *Assume that $\Sigma_{xx} > 0$ and $\gamma_1 > \cdots > \gamma_K > 0$ for some $K$.*

(i) *Under conditions $A_1$–$A_4$ in Theorem 1, the proposed filtering procedure will select all of the variables with nonzero coefficients in $\beta_i$ ($i = 1, 2, \ldots, K$) with probability tending to 1 as $n \to \infty$.*

(ii) *Suppose $\beta_i = (\beta_{i,1}, \beta_{i,2}, \ldots)$, where $\beta_{i,1}$ is a vector that contains the nonzero coefficients of $\beta_i$ and $\beta_{i,j} = 0$ ($j \geq 2$) are the remaining individual components. If condition $A_2$ is strengthened to the following $A_2'$, $A_2'$: $k_n^2/n \to 0$ and $k_n^4/n \to \infty$; then for $i = 1, \ldots, K$, $\sqrt{n}(\hat{\beta}_{i,1}^f - \beta_{i,1})$ has the same limiting distribution as $\sqrt{n}(\hat{\beta}_{i,1} - \beta_{i,1})$, and $\sqrt{n}|\hat{\beta}_{i,j}^f - \beta_{i,j}|$ is stochastically dominated by $\sqrt{n}|\hat{\beta}_{i,j} - \beta_{i,j}|$ for any $j \geq 2$.*

We refer to [8] for the limiting distribution of $\sqrt{n}(\hat{\beta}_i - \beta_i)$ under the unconstrained CANCOR method. In parametric settings, penalized likelihoods with $L_1$ constraints have been studied by earlier authors such as [6] and [7]. Here, a measure of correlation takes the role of a likelihood, and the nonparametric component (in the case of $k_n$ increasing with $n$) complicates the analysis. Indeed, the proposed dimension reduction is a combination of constrained canonical correlation ($C^3$), variable filtering, and re-estimation. The step in $C^3$ moves the estimated direction closer to $\beta_i$ under sparsity, which enables us to perform model selection through a simple variable filtering scheme. Re-estimation following variable filtering updates the direction estimates to achieve the highest possible correlation with a function of $y$. Theorem 2 shows that our proposed direction estimates are designed to explore the sparsity in $\beta_i$ with the same or better asymptotic efficiency than the CANCOR estimates.

The BIC-type criterion used in variable filtering is chosen to mimic the Bayesian information criterion, but it cannot be identified as the exact or approximate BIC in the usual sense. Rather, this choice is in part based on our empirical experience. The asymptotic results in Theorem 2 remain valid for a wide range of model selection criteria that balance correlation $r_d$ with complexity $d$.

**4. Simulation studies.** We perform 4 simulation studies in this section to study the performance of the proposed dimension reduction method, which consists of $C^3$, variable filtering and re-estimation. For convenience, in the simulation studies and in the examples in Sections 4 and 5, we will simply refer to the proposed dimension reduction as the $C^3$ method. To compare the $C^3$ method with the shrinkage sliced inverse regression (SSIR) method proposed in [17], we specify the settings of the first 3 studies to be the same as those in their paper.

TABLE 1
*Summary of Study 1*

| Method | $C^3$ | | SSIR | | | |
|---|---|---|---|---|---|---|
| Criterion | $\alpha = 0.01$ | $\alpha = 0.005$ | GCV | AIC | BIC | RIC |
| Sample size | | | $n = 60$ | | | |
| $A_3$ | 0.00 | 0.00 | 0.00 | 0.00 | 0.00 | 0.00 |
| (SE) | (0.00) | (0.00) | (0.00) | (0.00) | (0.00) | (0.00) |
| $A_{21}$ | 20.73 | 20.80 | 20.00 | 19.94 | 20.95 | 20.99 |
| (SE) | (0.06) | (0.05) | (0.11) | (0.11) | (0.02) | (0.01) |
| Sample size | | | $n = 120$ | | | |
| $A_3$ | 0.00 | 0.00 | 0.00 | 0.00 | 0.00 | 0.00 |
| (SE) | (0.00) | (0.00) | (0.00) | (0.00) | (0.00) | (0.00) |
| $A_{21}$ | 21.00 | 21.00 | 20.35 | 20.34 | 20.99 | 21.00 |
| (SE) | (0.00) | (0.00) | (0.08) | (0.08) | (0.01) | (0.00) |

As reviewed in Section 1, SSIR estimates a shrinkage index $\zeta = (\zeta_1, \ldots, \zeta_p)$ subject to $\|\zeta\|_{L_1} \leq t$. The central DRS $S_{y|x}$ is estimated by $\text{span}\{\text{diag}(\tilde{\zeta})\hat{B}\}$ in SSIR, where $\tilde{\zeta}$ is the estimate of $\zeta$ by SSIR and $\hat{B}$ is the matrix with columns as the e.d.r. direction estimates by SIR. The estimated shrinkage indices $\tilde{\zeta}$ compress some rows of $\hat{B}$ to 0, and the corresponding variables will not be involved in the dimension reduction results. There are four criteria used in SSIR for selecting the value of their tuning parameter $t$: the generalized cross validation (GCV) criterion, Akaike's information criterion (AIC), the Bayesian information criterion (BIC) and the residual information criterion (RIC). We refer to [17] for details on SSIR and those criteria.

In implementing $C^3$, we varied $\alpha$ in the set $\{0.025, 0.01, 0.005, 0.0025\}$ in the selection of the tuning parameter $t$, but found that the results were quite robust. For brevity, the results with $\alpha = 0.01$ and $0.005$ are reported here. To solve the optimization problems, we used an R interface with the FORTRAN subroutine VF13AD in the Harwell Subroutine Library. Other packages such as PROC NLP in SAS may be used for the same purpose.

In each study, we generate 100 data sets of the sample size 60 or 120. To generate the B-spline basis functions $\pi(y)$, the quadratic spline (order $m = 3$) with $k_n = 4$ internal knots is used. Accordingly, we use $I = 7$ slices in SSIR such that $m + k_n = I$. Since the explanatory variables are generated from normal distributions, the linearity condition holds in each study. The average numbers of 0 coefficients in the estimated constrained e.d.r. directions are summarized.

STUDY 1.

$$y = x_1 + x_2 + x_3 + 0.5\epsilon,$$

where $x = (x_1, \ldots, x_{24})^T \sim N(\underline{0}, I_{24})$, $\epsilon \sim N(0,1)$, and $x$ and $\epsilon$ are independent. In this study, the true direction is $\beta_1 = (1, 1, 1, 0, \ldots, 0)^T$ with 21 zero coefficients. Table 1 summarizes the average number of zero coefficients, as well as the corresponding standard errors, in $\text{diag}(\tilde{\alpha})\hat{\beta}_1$ estimated by SSIR and in $\hat{\beta}_1^f$ estimated by $C^3$, where:

- $A_3$ is the average number of zero coefficients out of the first 3 coefficients in $\text{diag}(\tilde{\alpha})\hat{\beta}_1$ for SSIR or in $\hat{\beta}_1^f$ for $C^3$;
- $A_{21}$ is the average number of zero coefficients out of the last 21 coefficients in $\text{diag}(\tilde{\alpha})\hat{\beta}_1$ or in $\hat{\beta}_1^f$.

Note that for the true direction $\beta_1$, we have $A_3 = 0$ and $A_{21} = 21$.

STUDY 2.

$$y = x_1 / \{0.5 + (x_2 + 1.5)^2\} + 0.2\epsilon,$$

where $x = (x_1, \ldots, x_{24})^T \sim N(\underline{0}, I_{24})$, $\epsilon \sim N(0,1)$, and $x$ and $\epsilon$ are independent. In this study, the true central DRS is the subspace spanned by $\beta_1 = (1, 0, \ldots, 0)^T$ and $\beta_2 = (0, 1, 0, \ldots, 0)^T$. Letting $\tilde{B} = \text{diag}(\tilde{\alpha})(\hat{\beta}_1, \hat{\beta}_2)$ estimated by SSIR and $\hat{B}^f = (\hat{\beta}_1^f, \hat{\beta}_2^f)$ estimated by $C^3$, Table 2 summarizes the averages $A_2$ and $A_{22}$, where:

- $A_2$ is the average number of zero rows out of the first 2 rows in $\tilde{B}$ for SSIR or in $\hat{B}^f$ for $C^3$;
- $A_{22}$ is the average number of zero rows out of the last 22 rows in $\tilde{B}$ or in $\hat{B}^f$.

Here, a zero row is a row vector with all elements equal to 0. Note that for the matrix $B = (\beta_1, \beta_2)$ corresponding to the true central DRS, $A_2 = 0$ and $A_{22} = 22$.

STUDY 3.

$$y = x_1 / \{0.5 + (x_2 + 1.5)^2\} + 0.2\epsilon,$$

where $\epsilon \sim N(0, 1)$, $x \in N(\underline{0}, \Sigma)$ with $\Sigma_{ij} = 0.5^{|i-j|}$ for $1 \leq i, j \leq 24$, and $\epsilon$ and $x$ are independent. The averages $A_2$ and $A_{22}$ defined in Study 2 are summarized in Table 3 with $A_2 = 0$ and $A_{22} = 22$ corresponding to the true central DRS.

TABLE 2
*Summary of Study 2*

| Method | $C^3$ | | SSIR | | | |
|---|---|---|---|---|---|---|
| Criterion | $\alpha = 0.01$ | $\alpha = 0.005$ | GCV | AIC | BIC | RIC |
| Sample size | | | $n = 60$ | | | |
| $A_2$ | 0.14 | 0.14 | 0.02 | 0.02 | 0.07 | 0.17 |
| (SE) | (0.03) | (0.03) | (0.01) | (0.01) | (0.03) | (0.04) |
| $A_{22}$ | 18.48 | 18.93 | 6.70 | 6.27 | 14.23 | 18.62 |
| (SE) | (0.23) | (0.17) | (0.25) | (0.25) | (0.26) | (0.16) |
| Sample size | | | $n = 120$ | | | |
| $A_2$ | 0.00 | 0.00 | 0.00 | 0.00 | 0.00 | 0.00 |
| (SE) | (0.00) | (0.00) | (0.00) | (0.00) | (0.00) | (0.00) |
| $A_{22}$ | 19.44 | 19.45 | 7.66 | 7.29 | 15.48 | 19.80 |
| (SE) | (0.14) | (0.15) | (0.27) | (0.27) | (0.23) | (0.14) |

STUDY 4.

$$y = x_1 + \cdots + x_{24} + 0.5\epsilon,$$

where $x = (x_1, \ldots, x_{24})^T \sim N(\underline{0}, I_{24})$, $\epsilon \sim N(0, 1)$, and $x$ and $\epsilon$ are independent. The true direction is $\beta_1 = (1, 1, \ldots, 1)^T$. Table 4 summarizes $A$ that is the average number of zero coefficients in $\operatorname{diag}(\tilde{\alpha})\hat{\beta}_1$ for SSIR or in $\hat{\beta}_1^f$ for $C^3$. Note that for the true direction $\beta_1$, we have $A = 0$.

In Studies 1–3, the residual information criterion (RIC) has the best overall performance among the four criteria of SSIR in identifying the relevant explanatory variables and filtering out the irrelevant variables. The proposed $C^3$ method works competitively to SSIR with RIC. However, in Study 4 where all of the variables are relevant, the $C^3$ method misses none of them while SSIR with RIC misses more than half of them on average. In Study 4, the AIC and GCV criteria in SSIR outperform BIC and RIC, although they still do not perform as well as the proposed $C^3$ method.

TABLE 3
*Summary of Study 3*

| Method | $C^3$ | | SSIR | | | |
|---|---|---|---|---|---|---|
| Criterion | $\alpha = 0.01$ | $\alpha = 0.005$ | GCV | AIC | BIC | RIC |
| Sample size | | | $n = 60$ | | | |
| $A_2$ | 0.27 | 0.33 | 0.07 | 0.07 | 0.18 | 0.22 |
| (SE) | (0.04) | (0.05) | (0.03) | (0.03) | (0.04) | (0.04) |
| $A_{22}$ | 18.94 | 19.37 | 5.79 | 5.28 | 14.35 | 18.52 |
| (SE) | (0.18) | (0.18) | (0.28) | (0.28) | (0.26) | (0.17) |
| Sample size | | | $n = 120$ | | | |
| $A_2$ | 0.18 | 0.25 | 0.01 | 0.01 | 0.03 | 0.06 |
| (SE) | (0.04) | (0.04) | (0.01) | (0.01) | (0.02) | (0.02) |
| $A_{22}$ | 20.10 | 20.22 | 6.77 | 6.29 | 14.65 | 18.62 |
| (SE) | (0.11) | (0.10) | (0.24) | (0.24) | (0.24) | (0.16) |

TABLE 4
*Summary of Study 4*

| Method | $C^3$ | | SSIR | | | |
|---|---|---|---|---|---|---|
| Criterion | $\alpha = 0.01$ | $\alpha = 0.005$ | GCV | AIC | BIC | RIC |
| Sample size | | | $n = 60$ | | | |
| A | 0.00 | 0.00 | 6.82 | 6.39 | 13.61 | 17.36 |
| (SE) | (0.00) | (0.00) | (0.21) | (0.21) | (0.18) | (0.17) |
| Sample size | | | $n = 120$ | | | |
| A | 0.00 | 0.00 | 0.15 | 0.09 | 5.91 | 13.81 |
| (SE) | (0.00) | (0.00) | (0.04) | (0.03) | (0.24) | (0.19) |

According to Shi and Tsai [18], the RIC criterion has a greater penalty function on the effective number of parameters than the AIC and BIC criteria. Thus, as shown in the simulation studies, SSIR with RIC usually compresses more coefficients to be 0 than SSIR with AIC or BIC. Therefore, in Studies 1–3 when only a small number of variables are involved in the central DRS, the RIC criterion helps filter out the irrelevant variables, but in Study 4 when all of the variables are involved, the RIC criterion is too aggressive, underscoring the difficulty in doing well for both types of situations. In $C^3$, instead of *just* penalizing the number of parameters, we monitor the decrease in the correlation, the objective function in CANCOR, which has a very intuitive appeal and makes it easy to outperform other criteria. In Study 4, excluding any of the variables will drive the correlation out of the lower limit of the confidence interval. As a result, $C^3$ correctly picks up all of the variables.

In summary, the selection criteria in SSIR have varied performance, but the proposed $C^3$ method competes favorably with the best criterion used in SSIR, and has a smaller risk for being overly aggressive. The simulation studies also show that the performance of both SSIR and $C^3$ can be improved as the sample size increases.

**5. Examples.** In this section, we apply the proposed $C^3$ method to two studies, the car price data and the Boston housing data.

5.1. *Example on car prices.* In the car price data analyzed in [16], the nonnegotiable transaction prices ($y$) of 25 different family saloons, together with their nine attributes, are recorded. Those nine attributes are mileage per gallon ($x_1$), horsepower ($x_2$), length ($x_3$), width ($x_4$), weight ($x_5$), height ($x_6$), satisfaction ($x_7$), reliability ($x_8$) and overall evaluation ($x_9$). Applying $C^3$ to the car price data, we aim to select a subset of the attributes to explain the variation among the transaction prices.

The nine attribute variables are standardized to mean 0 and variance 1. For estimating the dimensionality $K$ in model (1), the sequential tests are performed with different numbers of slices ($I$) in SIR and different numbers of internal knots ($k_n$) with quadratic spline (order $m = 3$) in CANCOR. To make the estimates $\hat{K}$ by SIR and CANCOR comparable, the numbers $I$ and $k_n$ are matched such that $I = m + k_n$. The dimensionality estimates by SIR and CANCOR with different $I$ and $k_n$, respectively, are summarized in Table 5. The estimate by SIR is sensitive to the number of slices while the estimate by CANCOR is robust against the number of internal knots. Thus, we select the CANCOR estimate $\hat{K} = 1$, and estimate the first constrained e.d.r. direction using $C^3$. For generating the B-spline basis functions based on $y$, the quadratic spline with $k_n = 4$ internal knots is used.

The first canonical correlation corresponding to the direction estimate $\hat{\beta}_1$ by CANCOR is 0.950 with the approximated 99.5% lower confidence limit 0.858. The iterative process described in Section 2 selects the tuning parameter $t = 1.10$. The corresponding constrained direction estimate at $t = 1.10$ yields the correlation 0.872, and has two nonzero coefficients selected by the variable filtering procedure. Using the two selected variables, the constrained direction is re-estimated as described in Section 3. The re-estimated constrained direction $\hat{\beta}_1^f$ yields the correlation 0.905. The direction estimates $\hat{\beta}_1$ and $\hat{\beta}_1^f$ by CANCOR and $C^3$, respectively, are standardized to unit length, and are shown in Table 6.

TABLE 5
*Dimensionality estimates for the car price data*

| SIR | | | | | | | CANCOR | | | | | | |
|---|---|---|---|---|---|---|---|---|---|---|---|---|---|
| $I$ | 3 | 4 | 5 | 6 | 7 | 8 | $k_n$ | 0 | 1 | 2 | 3 | 4 | 5 |
| $\hat{K}$ | 1 | 2 | 1 | 0 | 0 | 0 | $\hat{K}$ | 1 | 1 | 1 | 1 | 1 | 1 |

TABLE 6
*Direction estimates by CANCOR and $C^3$ for the car price data*

| | $x_1$ | $x_2$ | $x_3$ | $x_4$ | $x_5$ | $x_6$ | $x_7$ | $x_8$ | $x_9$ |
|---|---|---|---|---|---|---|---|---|---|
| $\hat{\beta}_1$ | 0.056 | 0.329 | $-0.472$ | 0.434 | $-0.032$ | 0.136 | 0.647 | $-0.142$ | $-0.138$ |
| $\hat{\beta}_1^f$ | 0 | 0.588 | 0 | 0 | 0 | 0 | 0.809 | 0 | 0 |

Back to the original scale of the variables, the estimated linear combination corresponding to $\hat{\beta}_1^f$ is $0.046x_2 + 0.999x_7$. The selected variables are horsepower ($x_2$) and satisfaction ($x_7$) by the proposed $C^3$ method. In [16], the SIR method is applied to the car price data, and the variables are ordered according to the direction estimate by SIR. Based on the ordering, a total number of 9 nested models are considered, and an AIC-based model selection criterion for single-index models is applied to the 9 nested models, resulting in the selected model $y = \hat{g}(0.045x_2 + 0.999x_7)$, where the function $\hat{g}(\cdot)$ is estimated by applying the local polynomial regression with the selected bandwidth. The direction estimated in [16] is almost identical to the one estimated by $C^3$. However, the $C^3$ method does not require an estimation of $g$, and it selects (rather than assumes) a uni-dimensional model that is adaptive to the data.

If we apply the variable filtering procedure to the CANCOR direction estimate $\hat{\beta}_1$, we end up with 4 variables with nonzero coefficients, $x_2$, $x_3$, $x_4$ and $x_7$. On the original scale of the variables, the corresponding estimated linear combination is $0.041x_2 - 0.171x_3 + 0.471x_4 + 0.864x_7$. Given that $x_3$ and $x_4$ have large coefficients in $\hat{\beta}_1$, they are forced in by the variable filtering procedure. However, when the $C^3$ method is performed, the coefficients of $x_3$ and $x_4$ decrease as the tuning parameter $t$ decreases, while the coefficient of $x_2$ increases until $t$ is very close to 1. When the iterative process stops at $t = 1.10$, the coefficients of $x_3$ and $x_4$ are very close to 0 in the constrained direction estimate. As a result, the variable filtering procedure does not select $x_3$ and $x_4$. This example confirms that when the explanatory variables are correlated, applying the simplified variable filtering procedure given in the paper to the unconstrained direction estimates would not be effective, and a fuller version of variable selection would be needed.

5.2. *Example with Boston housing data.* We now analyze the Boston housing data using the proposed $C^3$ method. Although this data set has been widely analyzed from different perspectives, our purpose is to demonstrate the difference between constrained and unconstrained direction estimates. The Boston housing data contains 506 observations, and can be downloaded from the web site 0http://lib.stat.cmu.edu/datasets/boston_corrected.txt. The dependent variable $y$ is the median value of owner-occupied homes in each of the 506 census tracts in the Boston Standard Metropolitan Statistical Areas. The 13 explanatory variables are per capita crime rate by town ($x_1$); proportion of residential land zoned for lots over 25,000 sq.ft ($x_2$); proportion of nonretail business acres per town ($x_3$); nitric oxides concentration ($x_4$); average number of rooms per dwelling ($x_5$); proportion of owner-occupied units built prior to 1940 ($x_6$); weighted distances to five Boston employment centers ($x_7$); full-value property-tax rate ($x_8$); pupil-teacher ratio by town ($x_9$); proportion of blacks by town ($x_{10}$); percentage of lower status of the population ($x_{11}$); Charles River dummy variable ($x_{12}$); index of accessibility to radial highways ($x_{13}$).

TABLE 7
*Direction estimates by CANCOR and $C^3$ for the Boston housing data*

|  | $\hat{\beta}_1$ | $\hat{\beta}_2$ | $\hat{\beta}_1^f$ | $\hat{\beta}_2^f$ |
|---|---|---|---|---|
| $x_1$ | 0.052 | −0.064 | 0 | 0 |
| $x_2$ | 0.085 | −0.266 | 0 | 0 |
| $x_3$ | −0.078 | −0.015 | 0 | 0 |
| $x_4$ | −0.009 | −0.352 | 0 | 0 |
| $x_5$ | 0.871 | −0.565 | 0.962 | −0.645 |
| $x_6$ | −0.306 | −0.058 | −0.174 | −0.096 |
| $x_7$ | −0.291 | −0.149 | 0 | 0 |
| $x_8$ | −0.165 | 0.022 | −0.166 | 0 |
| $x_9$ | −0.125 | −0.041 | −0.126 | 0 |
| $x_{10}$ | −0.005 | 0.089 | 0 | 0 |
| $x_{11}$ | 0.008 | −0.644 | 0 | −0.758 |
| $x_{12}$ | 0.043 | 0.108 | 0 | 0 |
| $x_{13}$ | 0.039 | 0.143 | 0 | 0 |

For observations with crime rate greater than 3.2, the variables $x_2$, $x_3$, $x_8$, $x_9$ and $x_{13}$ are constants except for 3 observations. Thus, as in [1] and [13], we use the 374 observations with crime rate smaller than 3.2 in this analysis. To make the explanatory variables comparable in scale, we standardize each of them individually to mean 0 and variance 1. The linearity condition is assumed given the arguments in [13]. As in the simulation studies and the car price data example, the quadratic spline (order $m = 3$) with $k_n = 4$ internal knots is used to generate the B-spline basis functions $\pi(y)$.

By CANCOR, there are four e.d.r. directions that are significant. A closer look at those directions showed that the 3rd and 4th directions are mainly due to outlying observations. To downweight the effect of outliers, we choose the weighted CANCOR approach of [21], where each observation is weighted according to its leverage in the $x$-space. With such weights, only two directions are found to be significant, and the resulting direction estimates by weighted CANCOR and the corresponding $C^3$ method, rescaled to unit length, are summarized in Table 7. Each of the estimated e.d.r. directions by the weighted CANCOR involves all of the explanatory variables. Due to the correlations among the variables, the coefficients cannot be interpreted individually. The directions estimated by $C^3$ are more focused, and 9 or 10 coefficients in each direction are compressed to 0.

Table 7 shows that, in the first constrained direction estimate $\hat{\beta}_1^f$, the variables $x_5$, $x_6$, $x_8$ and $x_9$ are singled out. However, the variable $x_5$ has the dominant loading. In the second constrained direction estimate $\hat{\beta}_2^f$, the variables $x_5$, $x_6$ and $x_{11}$ are singled out. Since the variable $x_5$, the housing size information, has already been picked up by $\hat{\beta}_1^f$, the second constrained e.d.r. direction $\hat{\beta}_2^f$ gets at how old and how wealthy the district is. We conclude that the housing size information and the district information (age and wealth) are well summarized by $x^T\hat{\beta}_1^f$ and $x^T\hat{\beta}_2^f$ to explain the variation of the median housing values among those 374 districts.

The direction estimates $\hat{\beta}_i$ by the weighted CANCOR and $\hat{\beta}_i^f$ by $C^3$ are quite close for the Boston housing data with correlations 0.982 between $x^T\hat{\beta}_1$ and $x^T\hat{\beta}_1^f$ and 0.849 between $x^T\hat{\beta}_2$ and $x^T\hat{\beta}_2^f$. However, with an acceptable sacrifice on the canonical correlations, the irrelevant or marginally important variables are filtered out, and the direction estimates by $C^3$ are easier to interpret.

**6. Conclusion.** The constrained canonical correlation ($C^3$) method is proposed to enhance the e.d.r. direction estimates by CANCOR. By imposing the $L_1$-norm constraint in the $C^3$ method, the noise contained in the data could be further filtered out and a small number of informative variables could be easily identified for each direction estimate. The directions estimated by the proposed method are easier to interpret and more helpful for the subsequent statistical analysis. Compared with the SSIR method and a recent work [10] on another shrinkage-type method, the proposed $C^3$ method has a unique advantage in that we use correlation as a transparent objective so that there is a natural criterion to guide us how much to "shrink." The proposed method results in final parameter estimates that are at least as efficient as the CANCOR estimates, but adapt nicely to sparsity in these directions whenever possible. Asymptotic properties for other shrinkage methods are unavailable for comparison at the moment.

An R interface with a Fortran subroutine is used to implement the proposed $C^3$ method in this paper. The programs can be obtained from the first author upon request.

## APPENDIX

PROOF OF THEOREM 1. Given a sample $\{X_t, Y_t\}_{t=1}^n$, we standardize $X_t$ to $Z_t = \hat{\Sigma}_{xx}^{-1/2}(X_t - \bar{X})$, where $(\bar{X}, \hat{\Sigma}_{xx})$ are the mean and covariance estimates. Using the variables $Z = (Z_1, \ldots, Z_n)^T$ as in [8], the direction estimates by CANCOR are $\hat{\beta}_i = \hat{\Sigma}_{xx}^{-1/2}\hat{\eta}_i$, where $\hat{\eta}_i$ are the eigenvectors of the matrix $\hat{\Delta}_n = n^{-1}Z^T\Pi^*(\Pi^{*T}\Pi^*)^{-1}\Pi^{*T}Z$ and $\Pi^*$ is the centered version of $\Pi = (\pi(Y_1), \ldots, \pi(Y_t))^T$. Letting $M = \text{Cov}(E(z|y)) = E[E(z|y)E(z^T|y)]$ and $(\eta_i, \lambda_i)$ are the eigenvectors and eigenvalues of $M$, we have $\beta_i = \Sigma_{xx}^{-1/2}\eta_i$ and $\gamma_i = \lambda_i^{1/2}$. It follows from the same arguments used in [8] that $\hat{\Delta}_n \xrightarrow{P} M$ and $\hat{\gamma}_i \xrightarrow{P} \gamma_i$.

To show $\hat{\beta}_i^c \xrightarrow{P} \beta_i$, it is equivalent to show that any subsequence $\{\hat{\beta}_{i,n_k}^c\}$ of $\{\hat{\beta}_i^c\}$, $n_k \subseteq \mathbb{N}$, contains a further subsequence that converges almost surely to $\beta_i$. We handle this for each $i$.

Let $(\hat{\alpha}_i, \hat{\beta}_i)$ and $(\hat{\alpha}_i^c, \hat{\beta}_i^c)$ be the direction estimates from CANCOR and $C^3$, respectively, and $\hat{\gamma}_i = \hat{\alpha}_i^T \hat{\Sigma}_{\pi x} \hat{\beta}_i$ and $\hat{\gamma}_i^c = \hat{\alpha}_i^{cT} \hat{\Sigma}_{\pi x} \hat{\beta}_i^c$. By the choice of the tuning parameter $t$, we have $(\hat{\gamma}_i^c - \hat{\gamma}_i) \to 0$, and thus $\hat{\gamma}_i^c \xrightarrow{P} \gamma_i$. Given the fixed dimensionality of $\Sigma_{xx}$, we also have $\hat{\Sigma}_{xx} \xrightarrow{P} \Sigma_{xx}$ under condition $A_3$. Therefore, for any subsequence $n_k$, we can find a further subsequence $n_{k'} \subseteq n_k$ such that $\hat{\gamma}_{i,n_{k'}}^c \to \gamma_i$ and $\hat{\Sigma}_{xx,n_{k'}} \to \Sigma_{xx}$ almost surely.

Let $O_i = \{\omega : \hat{\gamma}_{i,n_{k'}}^c(\omega) \to \gamma_i \text{ and } \hat{\Sigma}_{xx,n_{k'}}(\omega) \to \Sigma_{xx}\}$. It is clear that $P(O_i) = 1$. To prove Theorem 1, it suffices to show that $O_i \subseteq \{\omega : \hat{\beta}_{i,n_{k'}}^c(\omega) \to \beta_i\}$.

Otherwise, there exists some $\omega_0 \in O_i$ such that $\hat{\beta}_{i,n_{k'}}^c(\omega_0) \to \beta_i$ does not hold. Since $\Sigma_{xx} > 0$, $\hat{\Sigma}_{xx,n_{k'}}(\omega_0) \to \Sigma_{xx}$, and $\hat{\beta}_{i,n_{k'}}^{cT}(\omega_0)\hat{\Sigma}_{xx,n_{k'}}(\omega_0)\hat{\beta}_{i,n_{k'}}^c(\omega_0) = 1$, the sequence $\hat{\beta}_{i,n_{k'}}^c(\omega_0)$ must be bounded. Thus, there exists a subsequence $\hat{\beta}_{i,n_{k''}}^c(\omega_0) \to \beta^* \neq \beta_i$ with $n_{k''} \subseteq n_{k'}$ and $\beta^{*T}\Sigma_{xx}\beta^* = 1$.

Let $\Lambda(\omega_0) = \{\alpha : \alpha^T \hat{\Sigma}_{\pi\pi,n_{k''}}(\omega_0)\alpha = 1, \alpha^T \hat{\Sigma}_{\pi\pi,n_{k''}}(\omega_0)\hat{\alpha}_{l,n_{k''}}^c = 0, l = 1, \ldots, i-1\}$, and write $S_{n_{k''}} = \hat{\Sigma}_{\pi x,n_{k''}}(\omega_0)$ for simplicity, we have as $n_{k''} \to \infty$,

$$\max_{\alpha \in \Lambda(\omega_0)} |\alpha^T \hat{\Sigma}_{\pi x, n_{k''}}(\omega_0)(\hat{\beta}_{i,n_{k''}}^c(\omega_0) - \beta^*)| \to 0. \tag{2}$$

Furthermore, let

$$\alpha'_{n_{k''}}(\omega_0) = \arg\max_{\alpha \in \Lambda(\omega_0)} \alpha^T S_{n_{k''}} \beta^*,$$

$$\alpha''_{n_{k''}}(\omega_0) = \arg\max_{\alpha \in \Lambda(\omega_0)} \alpha^T S_{n_{k''}} \hat{\beta}_{i,n_{k''}}^c(\omega_0).$$

We then have

$$\alpha'^T_{n_{k''}}(\omega_0) S_{n_{k''}}(\hat{\beta}_{i,n_{k''}}^c(\omega_0) - \beta^*)$$
$$\leq \max_{\alpha \in \Lambda(\omega_0)} \{\alpha^T S_{n_{k''}} \hat{\beta}_{i,n_{k''}}^c(\omega_0)\} - \max_{\alpha \in \Lambda(\omega_0)} \{\alpha^T S_{n_{k''}} \beta^*\}$$
$$\leq \alpha''^T_{n_{k''}}(\omega_0) S_{n_{k''}} \hat{\beta}_{i,n_{k''}}^c(\omega_0) - \alpha''^T_{n_{k''}}(\omega_0) S_{n_{k''}} \beta^*$$
$$= \alpha''^T_{n_{k''}}(\omega_0) S_{n_{k''}} (\hat{\beta}_{i,n_{k''}}^c(\omega_0) - \beta^*),$$

which implies

$$\left| \max_{\alpha \in \Lambda(\omega_0)} \{\alpha^T S_{n_{k''}} \hat{\beta}_{i,n_{k''}}^c(\omega_0)\} - \max_{\alpha \in \Lambda(\omega_0)} \{\alpha^T S_{n_{k''}} \beta^*\} \right|$$
$$\leq \max_{\alpha \in \Lambda(\omega_0)} |\alpha^T S_{n_{k''}} (\hat{\beta}_{i,n_{k''}}^c(\omega_0) - \beta^*)|.$$

Since $\hat{\gamma}^c_{i,n_{k''}}(\omega_0) = \max_{\alpha \in \Lambda(\omega_0)} \{\alpha^T S_{n_{k''}} \hat{\beta}^c_{i,n_{k''}}(\omega_0)\}$, it follows from (2) that $\hat{\gamma}^c_{i,n_{k''}}(\omega_0) - \max_{\alpha \in \Lambda(\omega_0)} \{\alpha^T S_{n_{k''}} \beta^*\} \to 0$, and thus, $\alpha'^T_{n_{k''}}(\omega_0) S_{n_{k''}} \beta^* \to \gamma_i$.

Letting $h_{n_{k''}}(y) = \pi(y)^T_{n_{k''}} \alpha'_{n_{k''}}(\omega_0)$, where $\pi(y)_{n_{k''}}$ denotes the vector of the B-spline basis functions generated based on the $n_{k''}$ observations, we have $\widehat{\text{Cov}}(x^T \beta^*, h_{n_{k''}}(y)) \to \gamma_i$ as $n_{k''} \to \infty$, where $\widehat{\text{Cov}}$ is used in the proof to denote the sample covariance. We also have $\widehat{\text{Cov}}(h_{n_{k''}}(y), h_{n_{k''}}(y)) = 1$ by the construction of $h$, so

$$(3) \qquad \frac{\widehat{\text{Cov}}(x^T \beta^*, h_{n_{k''}}(y))}{[\widehat{\text{Cov}}(h_{n_{k''}}(y), h_{n_{k''}}(y))]^{1/2}} \to \gamma_l.$$

We rescale the continuous functions $h_{n_{k''}}(y)$ to get

$$g_{n_{k''}}(y) = \begin{cases} h_{n_{k''}}(y), & \text{if } \max_{y \in [a,b]} |h_{n_{k''}}(y)| \leq 1, \\ h_{n_{k''}}(y) \big/ \max_{y \in [a,b]} |h_{n_{k''}}(y)|, & \text{if } \max_{y \in [a,b]} |h_{n_{k''}}(y)| > 1, \end{cases}$$

to ensure $|g_{n_{k''}}(y)| \leq 1$ on $[a,b]$, and $0 < \widehat{\text{Cov}}(g_{n_{k''}}(y), g_{n_{k''}}(y)) \leq 1$. Replacing $h_{n_{k''}}(y)$ by $g_{n_{k''}}(y)$ in (3), we have

$$(4) \qquad \widehat{\text{Cov}}(x^T \beta^*, g_{n_{k''}}(y)) - \gamma_i [\widehat{\text{Cov}}(g_{n_{k''}}(y), g_{n_{k''}}(y))]^{1/2} \to 0.$$

Given condition $A_3$ and the fact $|g_{n_{k''}}(y)| \leq 1$, it is easy to see that the sample covariances converge to the population counterparts, and therefore $\text{Cov}(x^T \beta^*, g_{n_{k''}}(y)) - \gamma_i [\text{Cov}(g_{n_{k''}}(y), g_{n_{k''}}(y))]^{1/2} [\text{Cov}(x^T \beta^*, x^T \beta^*)]^{1/2} \xrightarrow{P} 0$, and

$$(5) \qquad \text{Corr}(x^T \beta^*, g_{n_{k''}}(y)) \xrightarrow{P} \gamma_i.$$

We define $R^2(\beta) = \max_{T(y) \in H} \text{Corr}^2(x^T \beta, T(y))$ for any $\beta \in \mathbb{R}^p$ with $\beta^T \times \Sigma_{xx} \beta = 1$, where $H$ is the set containing any transformation of $y$ with $E[T^2(y)] < \infty$. As in [1], it can be shown that $R^2(\beta) = \beta^T \text{Cov}[E(x|y)] \beta$. Since $\text{Cov}[E(x|y)]$ has distinct nonzero eigenvalues and the function $R^2(\beta)$ is continuous, we have that for any $\epsilon_1 > 0$, there exists $\epsilon_2 > 0$ such that

(6) $\quad \beta \in \Omega_i(\epsilon_1) \Rightarrow \text{Corr}(x^T \beta, T(y)) < \gamma_i - \epsilon_2 \quad$ for every $T(y) \in H$,

where $\Omega_i(\epsilon_1) = \{\beta : \beta^T \Sigma_{xx} \beta = 1, \beta^T \Sigma_{xx} \beta_l = 0, l \leq i-1, |\beta - \beta_i| > \epsilon_1\}$.

For $i = 1$, it is immediate to see that (5) contradicts to (6). We now consider $i > 1$ by iteratively using the fact that $\hat{\beta}^c_l \xrightarrow{P} \beta_l$ for $l \leq i-1$. By taking limits on both sides of $\hat{\beta}^{cT}_{i,n_{k''}}(\omega_0) S_{n_{k''}} \hat{\beta}^c_{l,n_{k''}}(\omega_0) = 0$, we get $\beta^* \Sigma_{xx} \beta_l = 0$ for $l \leq i-1$. Thus $\beta^*$ falls into $\Omega_i(\epsilon_1)$ for some positive $\epsilon_1$, which again leads to a contradiction to (6). $\square$

PROOF OF THEOREM 2. (i) If $\hat{\beta}'^c_i(d)$ (for some $d < p$) does not contain all the variables with nonzero coefficient in $\beta_i$, the consistency result in Theorem 1 and (6) imply that the corresponding $r_d$ will be strictly below $\gamma_i$ for sufficiently large $n$, and therefore the filtering procedure will favor the full set at $d = p$. Thus the result (i) of Theorem 2 follows.

(ii) As a result of (i), the union of the selected variables from the filtering procedure contains all of the variables involved in the e.d.r. directions $\{\beta_i\}$, with probability tending to 1.

Except for $\omega \in O_n$ with $P(O_n) = \delta_n \to 0$ as $n \to \infty$, the final direction estimate from the proposed method is the same as the CANCOR directions on the relevant variables (with nonzero coefficients in any of the $\beta_i$'s), that is, we have $\hat{\beta}^f_{i,1} = \hat{\beta}_{i,1}$ except for $\omega \in O_n$. Thus $|P\{\sqrt{n}(\hat{\beta}^f_{i,1} - \beta_{i,1}) \in C\} - P\{\sqrt{n}(\hat{\beta}_{i,1} - \beta_{i,1}) \in C\}| \leq \delta_n$ for any $C$, and $\sqrt{n}(\hat{\beta}^f_{i,1} - \beta_{i,1})$ has the same limiting distribution as that of $\sqrt{n}(\hat{\beta}_{i,1} - \beta_{i,1})$. On the zero coefficients $\beta_{i,j}$ with $j \geq 2$, it is either $\hat{\beta}^f_{i,j} = 0$ or $\hat{\beta}^f_{i,j} = \hat{\beta}_{i,j}$, so the stochastic dominance is immediate. $\square$

**Acknowledgments.** The authors thank an Associate Editor and the referees for their helpful comments and suggestions.

DEPARTMENT OF STATISTICS  
UNIVERSITY OF VIRGINIA  
CHARLOTTESVILLE, VIRGINIA 22904  
USA  
E-MAIL: jz9p@virginia.edu

DEPARTMENT OF STATISTICS  
UNIVERSITY OF ILLINOIS AT URBANA-CHAMPAIGN  
CHAMPAIGN, ILLINOIS 61820  
USA  
E-MAIL: x-he@uiuc.edu